\documentstyle[12pt]{article}
\begin{document}
\begin{titlepage}
\begin{flushright}
\end{flushright}
\vskip.3in

\begin{center}
{\Large \bf A Remark On  the FRTS realization  and Drinfeld Realization of Quantum
Affine Superalgebra $U_q(\hat { osp}(1,2))$} 
\vskip .2in

\large{Jintai Ding}
\vskip.2in
{\em Department of Mathematical Sciences,  University of Cincinnati}
\end{center}

\vskip 2cm
\begin{center}
{\bf Abstract}
\end{center}
In this paper, we present the hidden symmetry behind 
the  Faddeev-Reshetikhin-Takhtajan-Semenov-Tian-Shansky
realization  of quantum affine superalgebras $U_q(\hat { osp}(1,2)) $
and add  the q-Serre relation to the Drinfeld realization of 
$U_q(\hat { osp}(1,2)) $ \cite{GZ} derived from 
the FRTS realization. 


\end{titlepage}

%  Greek letters

\def\a{\alpha}
\def\b{\beta}
\def\d{\delta}
\def\e{\epsilon}
\def\g{\gamma}
\def\k{\kappa}
\def\l{\lambda}
\def\o{\omega}
\def\t{\theta}
\def\s{\sigma}
\def\D{\Delta}
\def\L{\Lambda}

\def\G{{ G}}
\def\Gk{{ G}^{(k)}}
\def\R{{ R}}
\def\hR{{\hat{ R}}}
\def\C{{\bf C}}
\def\P{{\bf P}}

\def\uqgh{{U_q[gl(m|n)^{(1)}]}}
\def\uqoh{{U_q[osp(1|2)^{(1)}]}}

% Shorthands for \begin{equation} and the like

\def\beq{\begin{equation}}
\def\eeq{\end{equation}}
\def\bea{\begin{eqnarray}}
\def\eea{\end{eqnarray}}
\def\ba{\begin{array}}
\def\ea{\end{array}}
\def\no{\nonumber}
\def\lt{\left}
\def\rt{\right}
\newcommand{\bq}{\begin{quote}}
\newcommand{\eq}{\end{quote}}

\newtheorem{Theorem}{Theorem}
\newtheorem{Definition}{Definition}
\newtheorem{Proposition}{Proposition}
\newtheorem{Lemma}[Theorem]{Lemma}
\newtheorem{Corollary}[Theorem]{Corollary}
\newcommand{\proof}[1]{{\bf Proof. }
        #1\begin{flushright}$\Box$\end{flushright}}

\newcommand{\sect}[1]{\setcounter{equation}{0}\section{#1}}
\renewcommand{\theequation}{\thesection.\arabic{equation}}

\sect{Introduction\label{intro}}

For quantum affine algebras, there are two types of loop 
realizations: 1) the FRTS realization and  2) the Drinfeld 
realization.  For the case of $U_q(\hat { gl}(n))$ \cite{DF}, 
the connection of these two realizations are established 
via the Gauss decomposition of L-operators. This method is 
recently used to derive Drinfeld realization for 
 quantum affine superalgebras. 
 \cite{GZ} \cite{Z}
\cite{CWWZ}. 

Though for  the case of  $U_q(\hat {osp}(1,2)) $, such a connection was established 
\cite{GZ}, there are two important aspects of the theory that are still needed to be addressed. 

The first aspect is the hidden symmetry of FRTS realization, namely in the FRTS realization,
there is a hidden symmetry of the generating L-operators implied in the hidden symmetry of 
the R-matrix used in the FRTS realization. This hidden symmetry  implies that the Drinfeld realization 
is a quotient of the FRTS realization.

Another aspect of the theory is that of 
q-Serre relation.   
Let $E(z)$ $F(z)$ and $H(z)$ be the generating current operators of 
the  affine Lie superalgebra $\hat { osp}(1,2)$. We know that 
the defining relations include the Serre relation: 
$$ E(z)(\{E(w), E(x)\})= (\{E(w), E(x)\})E(z), $$
 $$ F(z)(\{F(w), F(x)\})= (\{F(w), F(x)\})F(z), $$
which clearly is not implied by the relation: 
$$(z-w)E(z)E(w)=(w-z)E(w)E(z),$$
$$(z-w)F(z)F(w)=(w-z)F(w)F(z).$$
In this aspect, the Drinfeld realization of $U_q(\hat { osp}(1,2)) $ in 
\cite{GZ} is incomplete in the sense that the Serre relation is missing. 
With the help of the hidden symmetry, we derive a q-Serre relation for the 
generating current operators of 
$U_q(\hat { osp}(1,2))$.

This paper is organized as the following: in Section 2, we recall the main results 
about $U_q(\hat {osp}(1,2))$ in \cite{GZ}; and we present our main results in Section 3.

\section{The FRTS realization and the Drinfeld realization of 
$U_q(\hat { osp}(1,2))$.}

In this section, we will recall the main results and the notation in \cite{GZ}. 
The FRTS realization of  affine superalgebras starts with 
a super R-matrix,  $R(z)\in End(V\otimes V)$, where $V$ is a ${\bf Z}_2$ graded vector space,
and $R(z)$ satisfying the  condition
$R(z)_{\a\b,\a'\b'}\neq 0$  only when
$[\a']+[\b']+[\a]+[\b]=0$ mod$2$, and the
 ${\bf Z}_2$  graded Yang-Baxter equation (YBE) 
\beq\label{rrr}
R_{12}(z)R_{13}(zw)R_{23}(w)=R_{23}(w)R_{13}(zw)R_{12}(z).
\eeq
The multiplication for the tensor product is defined for
homogeneous elements $a,~ b,~ a'$, $b'$ by
\beq
(a\otimes b)(a'\otimes b')=(-1)^{[b][a']}\,(aa'\otimes bb'),
\eeq
where $[a]\in{\bf Z}_2$ 
denotes the grading of the element $a$. The FRTS realization of quantum affine 
superalgebras is given as the following.

\begin{Definition}\label{rs}: Let $R(\frac{z}{w})$ be a R-matrix
satisfying the   ${\bf Z}_2$   graded   YBE (\ref{rrr}). The FRTS superalgebra $U(\R)$ is
generated by invertible $L^\pm(z)$, satisfying
\bea
R({z\over w})L_1^\pm(z)L_2^\pm(w)&=&L_2^\pm(w)L_1^\pm(z)R({z\over
         w}),\no\\
R({z_+\over w_-})L_1^+(z)L_2^-(w)&=&L_2^-(w)L_1^+(z)R({z_-\over
         w_+}),\label{super-rs}
\eea
where $L_1^\pm(z)=L^\pm(z)\otimes 1$, $L_2^\pm(z)=1\otimes L^\pm(z)$
and $z_\pm=zq^{\pm {c\over 2}}$. For the first formula of
(\ref{super-rs}), the expansion direction of $R({z\over w})$ can be 
chosen in $z\over w$ or $w\over z$,  for the second formula, the
expansion direction is  in $z\over w$.
\end{Definition}

The algebra $U(\R)$ is a graded Hopf algebra: its coproduct is defined by
\beq
\D(L^\pm(z)=L^\pm(zq^{\pm 1\otimes {c\over 2}})\stackrel{.}{\otimes}
     L^\pm(zq^{\mp {c\over 2}\otimes 1}),
\eeq
and its antipode is
\beq
S(L^\pm(z))=L^\pm(z)^{-1}.
\eeq

For the case of $U_q(\hat { osp(1,2)})$, its FRTS realization is 
given with 
the super R-matrix 
$R({z\over w})\in End(V\otimes V)$, where 
$V$ is  the 3-dimensional vector representation of
$U_q({ osp(1,2)})$.
In V, we fix a set of basis vectors, $v_1,v_2$ and $v_3$, where
 $v_1,\;v_3$ are  graded 0 (mod 2) and $v_2$ is  graded 1(mod 2).
With this set of basis. the R-matrix is given as:
\beq
R({z\over w})=\left(
\begin{array}{ccccccccc}
1 & 0 & 0 & 0 & 0 & 0 & 0 & 0 & 0\\
0 & a & 0 & b & 0 & 0 & 0 & 0 & 0\\
0 & 0 & d & 0 & c & 0 & r & 0 & 0\\
0 & f & 0 & a & 0 & 0 & 0 & 0 & 0\\
0 & 0 & g & 0 & e & 0 & c & 0 & 0\\
0 & 0 & 0 & 0 & 0 & a & 0 & b & 0\\
0 & 0 & s & 0 & g & 0 & d & 0 & 0\\
0 & 0 & 0 & 0 & 0 & f & 0 & a & 0\\
0 & 0 & 0 & 0 & 0 & 0 & 0 & 0 & 1
\end{array}
\right), \label{r12}
\eeq
where 
\bea
&&a=\frac{q(z-w)}{zq^2-w},~~~~b=\frac{w(q^2-1)}{zq^2-w},~~~~
  c=\frac{q^{1/2}w(q^2-1)(z-w)}{(zq^2-w)(zq^3-w)},\no\\
&&d=\frac{q^2(z-w)(zq-w)}{(zq^2-w)(zq^3-w)},~~~~
  e=a-\frac{zw(q^2-1)(q^3-1)}{(zq^2-w)(zq^3-w)},\no\\
&&f=\frac{z(q^2-1)}{zq^2-w},~~~~
  g=-\frac{q^{5/2}z(q^2-1)(z-w)}{(zq^2-w)(zq^3-w)},\no\\
&&r=\frac{w(q^2-1)[q^3z+q(z-w)-w]}{(zq^2-w)(zq^3-w)},~~~~
  s=\frac{z(q^2-1)[q^3z+q^2(z-w)-w]}{(zq^2-w)(zq^3-w)}.
\eea
We also have that 
$R_{21}(\frac{z}{w})=R(\frac{w}{z})^{-1}$ . 

To derive the Drinfeld
current realization, we need the Gauss decomposition of 
the L-operators $L^\pm(z)$, which is given as 
\bea
L^\pm(z)&=&\left (
\begin{array}{ccc}
1 & 0 &  0\\
e^\pm_1(z) & 1 & 0\\
e^\pm_{3,1}(z) & e^\pm_2(z) & 1
\end{array}
\right )
\left (
\begin{array}{ccc}
k^\pm_1(z) & 0 & 0\\
0 & k^\pm_2(z) & 0\\
0 & 0 & k^\pm_3(z)
\end{array}
\right )
\left (
\begin{array}{ccc}
1 & f^\pm_1(z) & f^\pm_{1,3}(z) \\
0 & 1 &  f^\pm_2(z)\\
0 & 0 & 1
\end{array}
\right )
.\label{l+-}
\eea

Let   $\{X,Y\}\equiv XY+YX$ denotes an anti-commutator, and
\beq
\d(z)=\sum_{l\in {\bf Z}}\,z^l
\eeq
as a formal series.

Let  $X^\pm_i(z)$ be  defined as 
\bea
X^+_i(z)&=&f^+_{i,i+1}(z_+)-f^-_{i,i+1}(z_-),\no\\
X^-_i(z)&=&e^-_{i+1,i}(z_+)-e^+_{i+1,i}(z_-),
\eea
where $z_\pm=zq^{\pm{c\over 2}}$. 

Let 
\beq
X^\pm(z)=(q-q^{-1})\lt[X^\pm_1(z)+X^\pm_2(zq)\rt].\label{x=x1+x2}
\eeq

In \cite{GZ}, the following commutation relations are derived. 

\begin{Theorem} 

\bea
k^\pm_1(z)k^\pm_1(w)&=&k^\pm_1(w)k^\pm_1(z),\no\\
k^+_1(z)k^-_1(w)&=&k^-_1(w)k^+_1(z),\no\\
k^\pm_2(z)k^\pm_2(w)&=&k^\pm_2(w)k^\pm_2(z),\no\\
k^\pm_3(z)k^\pm_3(w)&=&k^\pm_3(w)k^\pm_3(z),\no\\
k^+_3(z)k^-_3(w)&=&k^-_3(w)k^+_3(z),\no\\
k^\pm_1(z)k^\pm_2(w)&=&k^\pm_2(w)k^\pm_1(z),\no\\
\frac{z_\pm-w_\mp}{z_\pm q^2-w_\mp} k^\pm_1(z)k^\mp_2(w)&=&
    \frac{z_\mp-w_\pm}{z_\mp q^2-w_\pm}
    k^\mp_2(w)k^\pm_1(z),\no\\
k^\pm_1(z)k^\pm_3(w)^{-1}&=&k^\pm_3(w)^{-1}k^\pm_1(z),\no\\
\frac{(z_\mp-w_\pm)(z_\mp q-w_\pm)}{(z_\mp q^2-w_\pm)(z_\mp q^3-w_\pm)}
   k^\pm_1(z)k^\mp_3(w)^{-1}&=&
\frac{(z_\pm-w_\mp)(z_\pm q-w_\mp)}{(z_\pm q^2-w_\mp)(z_\pm q^3-w_\mp)}
   k^\mp_3(w)^{-1}k^\pm_1(z),\no\\
\frac{z_\pm-w_\mp q}{z_\pm q-w_\mp}k^\pm_2(z)k^\mp_2(w)&=&
  \frac{z_\mp-w_\pm q}{z_\mp q-w_\pm}k^\mp_2(w)k^\pm_2(z),\no\\
k^\pm_2(z)^{-1}k^\pm_3(w)^{-1}&=&k^\pm_3(w)^{-1}k^\pm_2(z)^{-1},\no\\
\frac{z_\pm-w_\mp}{z_\pm q^2-w_\mp} k^\pm_2(z)^{-1}k^\mp_3(w)^{-1}&=&
    \frac{z_\mp-w_\pm}{z_\mp q^2-w_\pm}
    k^\mp_3(w)^{-1}k^\pm_2(z)^{-1},\label{k1k2k3}
\eea

\bea
k^\pm_1(z)X^-(w)k^\pm_1(z)^{-1}&=&\frac{z_\pm q^2-w}{q(z_\pm-w)}
    X^-(w),\no\\
k^\pm_1(z)^{-1}X^+(w)k^\pm_1(z)&=&\frac{z_\mp q^2-w}{q(z_\mp-w)}
    X^+(w),\no\\
k^\pm_2(z)X^-(w)k^\pm_2(z)^{-1}&=&\frac{(z_\pm-w q^2)(z_\pm q-w)}{q(z_\pm-w)
    (z_\pm-wq)}X^-(w),\no\\
k^\pm_2(z)^{-1}X^+(w)k^\pm_2(z)&=&\frac{(z_\mp-w q^2)(z_\mp q-w)}{q(z_\mp-w)
    (z_\mp-wq)}X^+(w),\no\\
k^\pm_3(z)X^-(w)k^\pm_3(z)^{-1}&=&\frac{z_\pm-wq^3}{q(z_\pm-wq)}
    X^-(w),\no\\
k^\pm_3(z)^{-1}X^+(w)k^\pm_3(z)&=&\frac{z_\mp-wq^3}{q(z_\mp-wq)}
    X^+(w),\label{x+-k1k2k3}
\eea
\bea
\frac{z-wq}{zq-w}X^-(z)X^-(w)+\frac{z-wq^2}{zq^2-w}X^-(w)X^-(z)
   &=&0,\no\\
\frac{z-wq^2}{zq^2-w}X^+(z)X^+(w)+\frac{z-wq}{zq-w}X^+(w)X^+(z)
   &=&0,\label{x++x--}
\eea
\bea
\{X^-(w), X^+(z)\}&=&\frac{-1}{q-q^{-1}}\lt[\d(\frac{z}{w}q^c)\lt(
   (1+q^{-{1\over 2}}-q^{1\over 2}) k^+_2(z_+)
    k^+_1(z_+)^{-1}-k^+_3(z_+q)k^+_2(z_+q)^{-1}\rt)\rt.\no\\
& &\lt.-\d(\frac{z}{w}q^{-c})\lt(k^-_2(w_+)k^-_1(w_+)^{-1}
   -(1+q^{-{1\over 2}}-q^{1\over 2})k^-_3(w_+q)k^-_2(w_+q)^{-1}\rt)\rt].\no\\
   \label{x+x-}
\eea
\end{Theorem}

In \cite{GZ},  the following definition is proposed: 
\bea
&&\phi_i(z)=k^+_{i+1}(z)k^+_i(z)^{-1},\no\\
&&\psi_i(z)=k^-_{i+1}(z)k^-_i(z)^{-1},~~i=1,2,\no\\
&&\phi(z)=(1+q^{-{1\over 2}}-q^{1\over 2})\phi_1(z)-\phi_2(zq),\no\\
&&\psi(z)=\psi_1(z)-(1+q^{-{1\over 2}}-q^{1\over 2})\psi_2(zq),\label{phi-psi}.
\eea
Then \cite{GZ},  
\begin{Theorem}:  $q^{\pm{c\over 2}},\;
X^\pm(z),\;\phi(z),\;\psi(z)$
give the defining relations of $U_q(\hat{osp}(1,2))$. More precisely,
$U_q(\hat {osp}(1,2))$ is an associative algebra with unit 1 and
the Drinfeld generators: $X^\pm(z),~\phi(z)$ and $\psi(z)$, a central
element $c$ and a nonzero complex parameter $q$. $\phi(z)$ and $\psi(z)$
are invertible. The gradings of the generators are: $[X^\pm(z)]=1$ and
$[\phi(z)]=[\psi(z)]=[c]=0$. The relations are given by
\begin{eqnarray*}
\phi(z)\phi(w)&=&\phi(w)\phi(z),\no\\
\psi(z)\psi(w)&=&\psi(w)\psi(z),\no\\
\phi(z)\psi(w)\phi(z)^{-1}\psi(w)^{-1}&=&\frac{(z_+q-w_-)(z_--w_+q)
    (z_+-w_-q^2)(z_-q^2-w_+)}{(z_+-w_-q)(z_-q-w_+)(z_+q^2-w_-)
    (z_--w_+q^2)},
\end{eqnarray*}
\begin{eqnarray*}
\phi(z)X^-(w)\phi(z)^{-1}&=&\frac{(z_+-w q^2)(z_+ q-w)}{(z_+ q^2-w)
    (z_+-wq)}X^-(w),\no\\
\phi(z)^{-1}X^+(w)\phi(z)&=&\frac{(z_--w q^2)(z_- q-w)}{(z_- q^2-w)
    (z_--wq)}X^+(w),\no\\
\psi(z)X^-(w)\psi(z)^{-1}&=&\frac{(z_--w q^2)(z_- q-w)}{(z_- q^2-w)
    (z_--wq)}X^-(w),\no\\
\psi(z)^{-1}X^+(w)\psi(z)&=&\frac{(z_+-w q^2)(z_+ q-w)}{(z_+ q^2-w)
    (z_+-wq)}X^+(w),\label{x-phipsi}
\end{eqnarray*}
\beq
\{X^+(z), X^-(w)\}=\frac{1}{q-q^{-1}}\lt[
     \d(\frac{w}{z}q^{c})\psi(w_+)
     -\d(\frac{w}{z}q^{-c})\phi(z_+)\rt].\label{x+x-=psiphi}
\eeq

\bea
\frac{z-wq}{zq-w}X^-(z)X^-(w)+\frac{z-wq^2}{zq^2-w}X^-(w)X^-(z)
   &=&0,\no\\
\frac{z-wq^2}{zq^2-w}X^+(z)X^+(w)+\frac{z-wq}{zq-w}X^+(w)X^+(z)
   &=&0.\label{Fx++x--}
\eea
\end{Theorem}

\section{The hidden symmetry of FRTS realization and the q-Serre relation of 
$U_q(\hat {osp}(1,2))$}

As we know, for the classical case,  there is a hidden symmetry on  
the three dimensional representation $V$ of 
$\hat {osp}(1,2)$, which comes from fact that 
V  is selfdual. 
It is , therfore, natural to expect that it is so for the quantum case. 

One observation on the results of \cite{GZ} in the previous section is that 
there are far too many generating current operators in the FRTS realization than 
that of  the Drinfeld realization. 
 This also 
clearly indicates that there should 
be a hidden symmetry of the L-operators, which can help us to resolve  this problem.

Let us start with the Heisenberg subalgebra of the FRTS realization. 
Through calculation (see also the formulas in the previous section), it is not difficult to derive
 that 
\begin{Proposition} 

$ k_1^{\pm}(z) k_3^{\pm}(zq^3)$ and 
$ k_2^{\pm}(z)( k_1^{\pm}(zq^{-1})(k_1^{\pm} (zq^{-2}))^{-1})^{-1} $ commute with 
$X^{\pm}_i(z)$. 

\end{Proposition} 

Clearly, from the point view of the theory of universal R-matrix \cite{D}\cite{FR}\cite{GZ}\cite{KT}, the image of  universal R-matrix of 
$U_q(\hat {osp}(1,2))$ on $V\otimes V$  actually is $f(z)R(z)$, where 
$f(z)$ is an analytic function. 
Therefore,  the $L^{\pm}(z)$ in the definition of $U(R)$ is already different by 
an multiple of an extra copy of Heisenberg algebra, which is the reason why 
in \cite{DF} we deal with $U_q(\hat  {gl}(n))$ not $U_q(\hat {sl}(n))$. 
Thus, 
we know that $U(R)$ is actually bigger. The difference comes from the 
Heisenberg subalgebra generated by $ k_1^{\pm}(z) k_3^{\pm}(zq^3)$. 

We can also derive the following 
\begin{Proposition}
$k_2^{\pm}(z)(k_1^{\pm}(z))^{-1}(
k_3^{\pm}(zq)
(k_2^{\pm}(zq))^{-1})^{-1}$ are central in $U(R)$. 
\end{Proposition}

From the point view of the universal R-matrix, we can see that 
these two central current operators are 
 nothing but identity operators. This shows that  
either $k_2^{\pm}(z)(k_1^{\pm}(z))^{-1}$ or $(
k_3^{\pm}(zq)
(k_2^{\pm}(zq))^{-1})$
as  generating operators for the subalgebra generated by 
$X^{\pm}_i(z)$. 

From now on, we   impose the condition that 

{\bf Condition I}
$$k_2^{\pm}(z)(k_1^{\pm}(z))^{-1}(
k_3^{\pm}(zq)
(k_2^{\pm}(zq))^{-1})^{-1}=1,$$
on the FRTS realization.

From this and by calculation, we can  deduce that 
 
\begin{Proposition} 
Let$$Y(z)=
X_2^{\pm}(z)-({\mp}q^{\mp \frac 1 2}X_1^{\pm}(zq^{-1}))=\Sigma_{n\in {\bf Z}} Y(n)z^{-n}. $$ 
Let $Y$ be the subalgebra  generated by  $Y(n), n\in {\bf Z}$.  
Then $YU(R)$ is an  double sided ideal in $U(R)$. 
\end{Proposition} 

 We can check that $YU(R)$ has no intersection with the subalgebra generated 
by $k^\pm_i(z)$, $z=,1,2,3$. 
If we look at the image of $Y(z)$ on V, we can see that it is actually zero.
If we look  from the point view of universal R-matrix\cite{KT},   we should have that 

{\bf  Condition II}
\bea Y(z)=0, \eea 
which, from now on,  we impose on the FRTS realization.

From the above two propositions, we can derive that 

\begin{Proposition} 
In the FRTS realization, 
$X^{\pm}_2(z)$, $k_2^+(z)(k_1^+(z))^{-1}$ and 
$k_2^+(z)(k_1^+(z))^{-1}$,  $ k_1^{\pm}(z) k_3^{\pm}(zq^3)$ and 
$ k_2^{\pm}(z)( k_1^{\pm}(zq^{-1})(k_1^{\pm} (zq^{-2}))^{-1})^{-1} $ 
  generate the whole algebra. 
\end{Proposition} 

This follows the similar calculation for the case $U_q(\hat { sl}(3))$\cite{DF} to  
 show $f^{\pm}_ {1,3}(z)$ and $e^{\pm}_{3,1}$ are generated by 
$X^{\pm}_2(z)$.

All the three propositions  and two conditions above gives us the hidden symmetry of the 
L-operators of the FRTS realization of $U_q(\hat {osp}(1,2))$. 
On the other hand, we should understand this hidden symmetry comes from the fact 
that on V, the representation of $U_q(\hat {osp}(1,2))$ is selfdual. Therefore
 $R(z)$ has a hidden symmetry like that of $U_q(o(n))$ and $U_q(sp(2n)$ in 
\cite{FRT} that basically  determines  the hidden symmetry of the L-operators\cite{KT}. 
This  should automatically lead us  to the  Condition (I) and (II).

As we explain in the introduction, from the classical theory, we know that 
the Drinfeld realization given in Theorem 2 in the section above is 
incomplete in the sense that q-Serre relation is missing. 
Before we deal with the q-Serre relation,  we would like to point out 
that the relation (2.14) and the similar ones in (2.18) are not completely 
right. The correct ones are 
given as following:

\begin{Proposition}

\bea
(zq^2-w)({z-wq})X^-(z)X^-(w)+(zq-w)({z-wq^2})X^-(w)X^-(z)
   &=&0,\no\\
(zq-w)({z-wq^2})X^+(z)X^+(w)+(zq^2-w)({z-wq})X^+(w)X^+(z)
   &=&0. \label{++}
\eea

\end {Proposition}

Here, the point is that the relation given in (2.14) and (2.18) are much 
stronger than (3.20) in the proposition above, and they actually imply the 
relations (3.20). However (3.20) do not imply 
(2.14) and (2.18).  The reason comes from the fact that they imply different 
pole conditions on the products $X^{\pm}(z)X^{\pm}(w)$. 
Similarly  the relations given in \cite{GZ} about $X_i^{\pm}(z)X_i^{\pm}(w)$ should 
be corrected correspondingly  as in the proposition above. 

With rather difficult calculation (similar to that in \cite{DF} for the case of 
$U_q(\hat { sl}(3))$), we can derive that 

\begin{Proposition}
$$ 
 \frac {(z_3-z_1q^{-1})(z_3-z_1q^{3})(z_1-z_2q^2)}{z_3-z_1q}X_i^+(z_3)X_i^+(z_1)X_i^+(z_2) +$$
$$
\frac { (z_2-z_1q^2) ({z_2-z_1q^{-1}}) (z_3-z_1q^{-1})(z_3-z_1q^{3})}{
(z_1-z_2q^{-1}) (z_3-z_1q) }X_i^+(z_3) X_i^+(z_2)X_i^+(z_1) )  -$$
$$
  \frac{((z_1-z_3q^2)(z_1-z_3q)(z_1q-z_3q^{-1})(z_1-z_2q^2)} {(z_1-z_3)(z_3-z_1q^2)}
X_i^+(z_1)X_i^+(z_2)X_i^+(z_3)- $$
\bea \frac{((z_1-z_3q^2)(z_1-z_3q)(z_1q-z_3q^{-1})(z_2-z_1q^2)(z_2-z_1q^{-1})}
 {(z_1-z_3)(z_3-z_1q^2)(z_1-z_2q^{-1})}
 X_i^+(z_2)X_i^+(z_1)X_i^+(z_3) )= 0 ,\no \\  \eea

$$ 
 \frac {(z_3-z_1q^{})(z_3-z_1q^{-3})(z_1-z_2q^{-2})}{z_3-z_1q^{-1}}X_i^-(z_3)X_i^-(z_1)X_i^-(z_2) + $$
$$\frac { (z_2-z_1q^{-2}) ({z_2-z_1q^{}}) (z_3-z_1q^{})(z_3-z_1q^{-3})}{
(z_1-z_2q^{}) (z_3-z_1q^{-1}) }X_i^-(z_3) X_i^-(z_2)X_i^-(z_1) )  -$$
$$
  \frac{((z_1-z_3q^{-2})(z_1-z_3q^{-1})(z_1q-z_3q^{})(z_1-z_2q^{-2})} {(z_1-z_3)(z_3-z_1q^{-2})}
X_i^-(z_1)X_i^-(z_2)X_i^-(z_3)- $$
\bea \frac{((z_1-z_3q^{-2})(z_1-z_3q^{-1})(z_1q-z_3q^{})(z_2-z_1q^{-2})(z_2-z_1q^{})}
 {(z_1-z_3)(z_3-z_1q^{-2})(z_1-z_2q^{})}
 X_i^-(z_2)X_i^-(z_1)X_i^-(z_3) )= 0,\no\\ \eea 
for $i=1,2,\emptyset$ and 
where the coefficient functions of the
relations above are expanded in the region of the expansion region 
of the corresponding monomial of the product of $X_i^{\pm}(z_j)$.
\end{Proposition} 

We call the two relations above, the q-Serre relations. 
It is not very difficult to show that this relation will degenerate into the 
classical Serre relations, but we still do not know how to write a simple 
formulation like that of the Drinfeld realization of $U_q(\hat { sl}(3))$.

\begin{Definition} $U_q(\hat {osp}(1,2))$ is an ${\bf Z}_2$ graded 
associative algebra 
generated by  $c$ , an central element; 
  $$\phi   (z)=\sum_{-m\in {\bf Z}_{\geq 0}}\phi   (-m)z^{m};$$ 
$$\psi   (z)=\sum_{m\in {\bf Z}_{\geq 0}}\psi   (-m)z^{m};$$
$$ X^\pm(z)=\Sigma X(n)z^{-n});$$ where 
 $\phi(z), \psi(z)$ are invertible and $\phi(0)\psi(0)=1=\psi(0)\phi(0).$
 The gradings of the generators are: $[\bar X^\pm(z)]=1$ and
$[\phi(z)]=[\psi(z)]=[c]=0$. The relations are given by (2.17),   (3.20), (3.22)and  (3.23).
\end{Definition}

\begin{Theorem} 
$U_q(\hat {osp}(1,2))$ is isomorphic to a  quotient of the subalgebra of the FRTS algebra generated by 
$X^{\pm}(z)$ and 
$k_3^\pm(z)k_2^\pm(z)^{-1}$, 
where the quotient ideal is generated by Condition (I) and (II), and the map is given by 
$$ \frac 1 {q-q^{-1}} X^{\pm}_2(zq) \rightarrow  X^\pm(z) ,$$
$$k_3^+(z)k_2^+(z)^{-1}\rightarrow \psi(z), $$ 
$$k_3^-(z)k_2^-(z)^{-1}\rightarrow \phi(z). $$
\end{Theorem} 

This theorem also tells us that the FRTS realization U(R) is nothing else  but 
$U_q(\hat{osp}(1,2))$ tensored by another copy of Heisenberg algebra. 

Starting from a completely different point of view,  we also come to the same definition of 
$U_q(\hat {osp}(1,2))$\cite{DFe}, where the q-Serre relation was derived in a different way.

\vskip.3in
\noindent {\bf Acknowledgments.} 
I would like to thank B. Feigin and S. Khoroshkin for useful discussions.

%\newpage

\end{document}